# ASYMPTOTICS FOR HITTING TIMES

By M. Kupsa and Y. Lacroix

*Charles University and Université de Toulon et du Var*

In this paper we characterize possible asymptotics for hitting times in aperiodic ergodic dynamical systems: asymptotics are proved to be the distribution functions of subprobability measures on the line belonging to the functional class

(A) $\quad \mathcal{F} = \left\{ F : \mathbb{R} \to [0,1] : \begin{bmatrix} F \text{ is increasing, null on } ]-\infty, 0]; \\ F \text{ is continuous and concave;} \\ F(t) \leq t \text{ for } t \geq 0. \end{bmatrix} \right\}.$

Note that all possible asymptotics are absolutely continuous.

**1. Introduction.** Throughout $(X, \mathcal{B}, \mu)$ is a *probability space*, $T : X \to X$ is measurable and preserves $\mu$. We also assume the *dynamical system* $(X, T, \mu)$ to be *ergodic*.

For $U \subset X$ with $\mu(U) > 0$, *Poincaré's recurrence theorem* states that the variable

$$\tau_U(x) = \inf\{k \geq 1 : T^k x \in U\}$$

is $\mu$-a.s. well defined. If $x \in U$, $\tau_U(x)$ denotes the *return time* of $x$ to $U$, and for arbitrary $x \in X$, $\tau_U(x)$ is the *hitting time* of $x$ to $U$ (also often called entrance time). *The return time theorem* [Kac (1947)] reads

$$\mathbb{E}(\mu(U)\tau_U) = \sum_{t \geq 1} t\mu(U \cap \{\tau_u = t\}) = 1,$$

where the expectation is computed with respect to the induced probability measure on $U$, $\mu_U := \frac{\mu}{\mu(U)}$.

Finer statistical properties of the variable $\mu(U)\tau_U$ have been investigated; for instance, Chazottes (2003) states conditions for the existence of higher-order moments, in connection with mixing properties of the system.

Another approach, rapidly developing in the last decade, relevant to the study of *recurrence to rare events in dynamical systems*, is to describe asymptotics for hitting or return times.

We say a sequence of distribution functions $(F_n)$ converges weakly to a function $F$ (which might not be a distribution function itself) if $F$ is increasing (not necessarily strictly) and at any point of continuity of $F$, say $t_0$, $F_n(t_0) \to F(t_0)$. Notice that we assume $F$ increasing a priori. We will write $F_n \Rightarrow F$ if $(F_n)$ converges weakly to $F$.

For $U \subset X$ measurable with $\mu(U) > 0$, we define

$$\widetilde{F}_U(t) := \frac{1}{\mu(U)} \mu(U \cap \{\tau_U \mu(U) \le t\})$$

and

$$F_U(t) = \mu(\{\mu(U)\tau_U \le t\}).$$

Let $(U_n)_{n \ge 1}$ denote a sequence in $\mathcal{B}$ with $\mu(U_n) \to 0$. The question of asymptotics asks for weak convergence of $(\widetilde{F}_{U_n})_{n \ge 1}$ or $(F_{U_n})_{n \ge 1}$, and in the case it does, asks for the nature of the limit. The latter concerns hitting times, and the former return times.

Weak limits (for both hitting and return times) have been shown to exist for suitably chosen $(U_n)_{n \ge 1}$ (essentially decreasing sequences of balls in a metric space $X$) and identified to be the distribution function of the positive exponential law with parameter 1, in many classes of mixing systems, in Abadi and Galves (2001), Collet and Galves (1993), Hirata, Saussol and Vaienti (1999), Saussol (1998) and Young (1999). Nonexponential asymptotics have been obtained in Coelho and de Faria (1990), for instance. The literature on the subject is rather important and our list is incomplete. For further information we refer the reader to the surveys Abadi (2004) or Coelho (2000).

There exists a connection between return time asymptotics and hitting time asymptotics, indeed, as shown in Hirata, Saussol and Vaienti (1999), when the asymptotics for return times is exponential with parameter 1, then so is the one for hitting times.

Possible asymptotics for return times, that is, weak limits for $(\widetilde{F}_{U_n})_{n \ge 1}$, were determined in Lacroix (2002). These are $\widetilde{F}$'s in $[0,1]$, null on $]-\infty, 0]$, increasing, such that $\int_0^{+\infty}(1 - \widetilde{F}(t))\,dt \le 1$.

Though asymptotics for hitting times have been studied in many papers in the literature, the question of the nature of possible asymptotics is still completely open.

We answer this question from probability theory. Let

(A) $\quad \mathcal{F} = \left\{ F : \mathbb{R} \to [0,1] : \begin{array}{l} F \text{ is increasing, null on } ]-\infty, 0]; \\ F \text{ is continuous, concave on } [0, +\infty[; \\ F(t) \le t \text{ for } t \ge 0. \end{array} \right\}.$

Notice that $\mathcal{F}$ contains only absolutely continuous distributions, some of which are associated to subprobability measures.

We say the system $(X, T, \mu)$ is *aperiodic* if for any $m \geq 1$,
$$\mu(\{x : T^m x = x\}) = 0.$$
We prove:

THEOREM 1. *Let $(X, T, \mu)$ be an ergodic aperiodic dynamical system. Given $F : \mathbb{R} \to \mathbb{R}$ increasing, there exists $(U_n)_{n \geq 0}$ with $\mu(U_n) \to 0$ and $F_{U_n} \Rightarrow F$ if and only if $F \in \mathcal{F}$.*

*Hence possible asymptotics for hitting times are exactly the elements of $\mathcal{F}$.*

We stress that the class $\mathcal{F}$ is rather restricted, which is unexpected.

Let us remark that the continuous parameter case has been studied in this journal, namely in Geman (1973). The characteristics of possible asymptotics in that case differ from ours, in that, for instance, asymptotics for the continuous parameter case may have a discontinuity jump at the origin. The proof technique is quite different, too.

A SHORT SKETCH OF THE (ELEMENTARY) PROOF OF THEOREM 1. Our proof uses the same techniques as those developed in Lacroix (2002). We provide, however, a few simplifications. We think that once the spirit of the proof is understood, details are easy to follow.

Here is how the proof goes: first we state (conditions $\mathcal{C}$ in Section 2) necessary conditions for an $F$ to be an $F_U$ for some $U$ in some ergodic system. We then define rational $F$'s, which are those satisfying $\mathcal{C}$ with additional rationality assumptions. These rational $F$'s are shown in the stamp machine lemma to be exactly those arising from periodic ergodic systems.

Second, the concavity–continuity lemma characterizes weak limits of rational $F$'s as to be exactly the elements of the class $\mathcal{F}$ described above.

Third, in a periodic system a $U$ is a finite collection of points, spaced along the irreducible cycle that defines the ergodic periodic transformation. It defines spacing and return times, and the rational $F$'s thereby produce models (that we call stamps) that enable one to mark the levels in Rohlin towers [cf. Shields (1973) for definitions]. Then if we call $\widetilde{U}$ the union of the marked levels in the tower, it is easy to see that the larger the tower is, the uniformly closer $F_U$ and $F_{\widetilde{U}}$ are. This is the approximation lemma, using the classical Rohlin lemma.

Fourth, there is also an obvious observation that any $F_U$ is arbitrarily uniformly close to a rational $F$. It follows at once that asymptotics in a given aperiodic ergodic system are the same as weak limits of $F_U$'s arising from all periodic ergodic systems.

Let us proceed.

**2. Kac's towns, stamps and Rohlin towers.** Since $\mu(U) > 0$, it is a standard construction in ergodic theory to build Kac's town above $U$, which is a juxtaposition of skyscrapers: the ground $U$ is partitioned into sets $U \cap \{\tau_U = k\}$, $k \geq 1$, and above each of these the action of $T$ goes upward along the floors of a skyscraper of height $k-1$. Once reached, the top floor points go back somewhere in $U$ under the action of $T$.

The measure of the union of the floors in Kac's town, including the ground floor, equals 1; this is yet another expression of Kac's return time theorem.

The connection with hitting time can be made as follows: if one wants to compute $\mu(\tau_U = k)$ for some $k \geq 1$, one has to take the measure of the union of the floors that have $k-1$ levels left above.

A closer look at an $F_U$ shows that it must have the following elementary properties:

1. Its discontinuities are located at points $\mu(U), 2\mu(U), \ldots$, the collection of which might be finite or not, depending on the fact that entry times to $U$ are bounded or not.
2. The distribution function $F_U$ is simple (it is the distribution function of a discrete random variable), constant on intervals of length $\mu(U)$ [the random variable concerned is $\mu(U)\tau_U$], is 0 on $]-\infty, \mu(U)[$ and tends to 1 at $+\infty$.
3. It has decreasing jumps of discontinuity; this is because the value of the jump at point $k\mu(U)$ equals the measure of the union of the floors having $k-1$ floors left above in Kac's town, which necessarily decreases with $k$.
4. The first discontinuity jump equals $\mu(U)$.

The conditions enumerated above—denote them by the symbol $\mathcal{C}$—are necessary for a distribution function to be an $F_U$ for some $U$ of positive measure and for some ergodic system $(X, T, \mu)$.

We will need the following definition: a distribution function $F$ on the real line is rational if it satisfies $\mathcal{C}$, has finitely many discontinuity points, all located at rationals, and has rational discontinuity jumps.

STAMP MACHINE LEMMA. *A distribution function $F$ is an $F_U$ for some ergodic periodic system $(X, T, \mu)$ if and only if it is rational.*

PROOF. The necessity follows from the preceding discussion, and the fact that in a periodic ergodic system, the set of possible return times is finite, any $U$ with positive measure has rational measure, and for any such $U$ any floor in Kac's town has a rational measure.

Conversely, given a rational $F$, we will build Kac's town with base $U$ such that $F = F_U$. The first thing to do is to collect the collection of decreasing discontinuity jumps in decreasing order, and to sort out from this collection the set of (decreasing) values of jumps, and for each value of a jump, the cardinality of the consecutive run of jumps having the given selected value.

The first set will provide the opportunity to compute the measures of the floors of the skyscrapers, while the second one will provide the possibility to compute the heights of the skyscrapers.

Let us assume that the discontinuities of $F$ are located at points $\alpha, 2\alpha, \ldots, K\alpha$, with $\alpha \in \mathbb{Q}$. We have a discontinuity jump $\beta_k = F(k\alpha^+) - F(k\alpha^-) \in \mathbb{Q}$, for $1 \leq k \leq K$. Notice that $\beta_1 = \alpha$, and that

$$\beta_1 + \cdots + \beta_K = 1,$$

because $F$ goes from 0 to 1 upward along its discontinuity jumps.

There are some $s \geq 1$ and integers $1 \leq k_1 < k_2 < \cdots < k_s = K$ such that

$$(\beta_1, \ldots, \beta_K) = (\beta_1 = \cdots = \beta_{k_1} > \beta_{k_1+1} = \cdots = \beta_{k_2} > \cdots > \beta_{k_{s-1}+1} = \cdots = \beta_{k_s}).$$

We draw a ground floor of measure $\alpha$, as well as, piled vertically underneath, downward, $k_1 - 1$ underground floors, of the same measure $\alpha = \beta_{k_1}$.

Next we pile rightmost, underneath, downward, $k_2$ consecutive floors of measure $\beta_{k_2}$, next $k_3$ ones, the same way, and so on.

The procedure ends up with something looking like Kac's town, but mirrored downward. Never mind; we reverse direction, and get the picture of something looking much like Kac's town. The union of the floors has measure 1 since $\beta_1 + \cdots + \beta_K = 1$. It now remains to find a periodic ergodic system, together with a $U$, for which this construction is the construction of Kac's town associated to $U$. We can write $\alpha = p/q$ and $\beta_{k_j} = p_j/q$, $1 \leq j \leq s$, for some denominator $q$.

We denote by $(X, T, \mu)$ the periodic ergodic system with $q$ elements and period $q$. We will construct $U \subset X$ with $p$ elements: we set $X = \{1, 2, 3, \ldots, q\}$, and $Tx = x + 1$ if $x < q$, $Tq = 1$. The measure $\mu$ is the equidistribution. We set

$$\begin{aligned}
U = \{&1, 1 + k_1, \ldots, 1 + (p_1 - p_2)k_1, \\
&1 + (p_1 - p_2)k_1 + k_2, \ldots, 1 + (p_1 - p_2)k_1 + (p_2 - p_3)k_2, \\
&\ldots, \ldots, \ldots, \ldots, \ldots, \ldots, \ldots, \\
&1 + (p_1 - p_2)k_1 + \cdots + (p_{s-1} - p_s)k_{s-1}, \\
&1 + (p_1 - p_2)k_1 + \cdots + (p_{s-1} - p_s)k_{s-1} + k_s, \ldots, \\
&1 + (p_1 - p_2)k_1 + \cdots + (p_{s-1} - p_s)k_{s-1} + p_s k_s\}.
\end{aligned}$$

Then $U$ contains $p = p_1 = 1 + (p_1 - p_2) + \cdots + (p_{s-1} - p_s) + p_s - 1$ elements, hence $\mu(U) = \alpha$.

Also notice that since

$$\begin{aligned}
q &= q(\beta_1 + \cdots + \beta_K) \\
&= (k_1(p_1 - p_2) + \cdots + k_{s-1}(p_{s-1} - p_s) + k_s(p_s - 1) + k_s),
\end{aligned}$$

possible return times to $U$ are exactly $k_1, \ldots, k_s$, and exactly $p_1 - p_2$ elements of $U$ return to $U$ at time $k_1$, $p_2 - p_3$ of them do so for time $k_2$, and so on.

We think that now the best thing to convince the reader that such system with such $U$ makes $F = F_U$ is to let him work out a handmade example along with the above guidelines, maybe also using the example developed below. □

DEFINITION 1. Given rational $F$, with parameters $K$, $\alpha = p/q$, $\beta_{k_j} = p_j/q$, we can construct $(X, T, \mu)$ periodic and $U$ as in the proof of the preceding lemma.

A stamp for $F$ marks the levels in a tower of a given width and a height equal to $q$ by marking those that have the same heights (counting from the stamps base) as those of the floors corresponding to heights in $U$ in the preceding lemma.

AN EXAMPLE OF A STAMP CONSTRUCTION. With the above notation let $F$ be rational with parameters $\alpha = 5/27$ and $(\beta_1 \geq \beta_2 \geq \cdots \geq \beta_7) = (5/27, 5/27, 5/27, 3/27, 3/27, 3/27, 3/27)$.

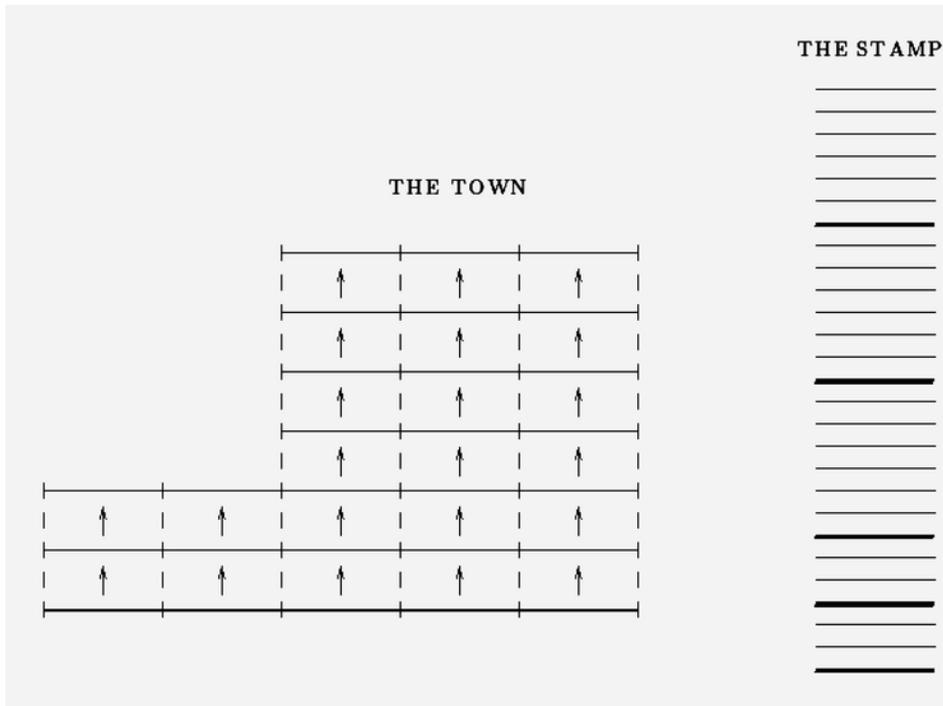

FIG. 1. *From $F$ to the town to its stamp.*

In Figure 1, for a finite periodic system with 27 points, we show Kac's town with base $U$ consisting of five points, where each floor in the town consists of a single point in the space. Each level in the town corresponds to a measure equal to the corresponding $\beta_i$, which reads, for our example, successively $5, 5, 5, 3, 3, 3, 3$ points. We also figure how this town produces stamps for $F$: we mark the base $U$, and unmark the floors above, then we pile the skyscrapers right above left from left to right. We obtain a vertical tower, in which marked floors (thickened in Figure 1) have a union equal to $U$.

**3. Weak limits and the class $\mathcal{F}$.** The proof of the following lemma is left to the reader; it may look very much like the one of the last statement of the concavity–continuity lemma below.

WEAK LIMITS LEMMA. *Given any $\varepsilon > 0$, any ergodic aperiodic system $(X, T, \mu)$ and any $U \subset X$ with $\mu(U) > 0$, there exists a rational $F$ such that for any $t \geq 0$, there exists an $s \geq 0$ with*

$$|s - t| < \varepsilon \quad \text{and} \quad |F(s) - F_U(t)| \leq \varepsilon.$$

Next we have:

CONCAVITY–CONTINUITY LEMMA. *Any weak limit of a sequence $(F_{U_n})$ arising from a system $(X, T, \mu)$ and some sequence of sets of positive measures in it, $(U_n)$, with $\mu(U_n) \to 0$, must be in $\mathcal{F}$ [cf. (A)].*
*Any $F \in \mathcal{F}$ is the weak limit of a sequence of rational $F$'s.*

PROOF. Assume $F_{U_n} \Rightarrow F$. Then $F$ must be increasing, hence has an at most countable set of discontinuity points. It has a dense set of continuity points whence $F$ must be zero on $]-\infty, 0[$ and must take its values in $[0, 1]$.

Recall that $F_{U_n}$ satisfies conditions $\mathcal{C}$. In particular, if $0 \leq s < t$, the increase $F_{U_n}(t) - F_{U_n}(s)$ is at most equal to the maximal discontinuity jump of $F_{U_n}$, that is, $\mu(U_n)$, times the number of intervals of length $\mu(U_n)$ needed to go from $s$ to $t$, plus 1. This means we have an inequality

(I) $\qquad 0 \leq F_{U_n}(t) - F_{U_n}(s) \leq \mu(U_n)\left(\frac{t-s}{\mu(U_n)} + 1\right) = t - s + \mu(U_n).$

We have $\mu(U_n) \to 0$, and the continuity modulus of $F_{U_n}$ is

$$\delta(F_{U_n}) := \limsup_{\varepsilon \downarrow 0^+} \sup_{|x-y| < \varepsilon} |F_{U_n}(x) - F_{U_n}(y)| = \mu(U_n).$$

Let us suppose $F$ has a discontinuity point at $x_0 \geq 0$; because $F$ increases, there exists $F(x_0^-) = \lim_{y < x_0} \uparrow F(y) < F(x_0^+) := \lim_{y > x_0} \downarrow F(y)$. Let us put $\delta = \delta(F, x_0) := F(x_0^+) - F(x_0^-)$. There exists $x_0 - \frac{\delta}{4} < x_1 < x_0 < x_2 < x_0 + \frac{\delta}{4}$ such that $F$ is continuous at $x_1$ and $x_2$. Then $F_{U_n}(x_1) \to F(x_1)$ and $F_{U_n}(x_2) \to F(x_2)$.

Passing to the limit and applying inequality (I), we obtain

$$0 < \delta \leq F(x_2) - F(x_1) = \lim_n (F_{U_n}(x_2) - F_{U_n}(x_1))$$

$$\leq \limsup_n (x_2 - x_1) + \mu(U_n) = x_2 - x_1 < \frac{\delta}{2},$$

a contradiction. So $F$ is continuous.

The fact that the weak limit is continuous makes it a simple limit. Each $F_{U_n}$ clearly satisfies $F_{U_n}(t) \leq t$ for $t \geq 0$ (cf. conditions $\mathcal{C}$). So the simple limit $F$ satisfies $F(t) \leq t$ for $t \geq 0$, too.

The concavity of $F$ is a consequence of the fact that $F_{U_n}$'s have decreasing discontinuity jumps.

To prove that an $F_0 \in \mathcal{F}$ is the weak limit of a sequence of rational $F$'s, it is enough, because $F_0$ is continuous, to prove that for any $\varepsilon > 0$, there exists a rational $F$ such that for any $\frac{1}{\varepsilon} \geq t \geq 0$, there exists an $s \geq 0$ with

$$(\star) \qquad |s - t| < \varepsilon \quad \text{and} \quad |F_0(t) - F(s)| \leq \varepsilon.$$

This can be done about the same way we would be proving the weak limits lemma, except we have a more complicated situation with a truncation, because $F_0$ might not grow up to 1 at $+\infty$.

Never mind; pick an integer $N > \frac{1}{\varepsilon}$, and divide the interval $[0, N]$ into intervals of equal lengths $\frac{1}{N}$. This produces a sequence of decreasing jumps $F_0(\frac{k+1}{N}) - F_0(\frac{k}{N})$, $0 \leq k < N$ (by concavity of $F_0$). We then approximate each $F_0(\frac{k}{N})$, $k > 0$, from above, by a positive rational at distance at most $\varepsilon$, and less than 1, call it $F(\frac{k}{N})$, in such a way that the new sequence of jumps $F(\frac{k+1}{N}) - F(\frac{k}{N})$ still decreases, with $F(0) = 0$.

Then we put $F(s) = F(\frac{[Ns]}{N})$. Finally, we complete if necessary by growing up from $F(N)$ to 1 using smaller rational jumps, being still constant on intervals of length $\frac{1}{N}$.

The obtained $F$ is rational and satisfies $(\star)$. $\square$

## 4. Stamping along Rohlin towers.

APPROXIMATION LEMMA. *Given $(X, T, \mu)$ ergodic and aperiodic, given a rational $F$ and given $\varepsilon > 0$, there exists $U \subset X$ measurable with $\mu(U) > 0$ and such that for any $t \geq 0$, there exists an $s \geq 0$ with both $|s - t| \leq \varepsilon$ and $|F(t) - F_U(s)| \leq \varepsilon$.*

PROOF. The rational $F$ may be realized in a periodic system of period $q > 0$, and produces stamps of arbitrary widths and height $q$, by the stamp machine lemma. We assume $F$ comes along with its rational parameters $\alpha, \beta_1 \geq \cdots \geq \beta_K$, and $s \geq 1$ such that

$$\beta_1 = \cdots = \beta_{k_1} > \beta_{k_1+1} = \cdots > \beta_{k_{s-1}+1} = \cdots = \beta_{k_s}$$

(we take the notation from Section 2). We denote $\alpha = p/q$, $\beta_{k_j} = \frac{p_j}{q}$, $1 \leq j \leq s$.

Our construction of a stamp for this $F$ in fact produces a subset $W$ consisting of $p$ points in a periodic system of period $q$, all spaced one with the following one by some $k_j$, such that $p_j - p_{j+1}$ such "spacings" equal $k_j$ for each $1 \leq j \leq s$ ($p_{s+1} = 0$). We assume the first point belongs to $W$ (as in Section 2).

By ergodicity and aperiodicity, using the Rohlin lemma, given $\delta > 0$ and $\frac{1}{\delta} \leq \frac{N}{q}$, there exist $V \subset X$ and $n \geq N$ with $0 < q\mu(V) < \delta$, $V, TV, \ldots, T^{n-1}$ disjoint, and $\mu(\bigcup_{k<n} T^k V) \geq 1 - \delta$.

We can write $n = qr + s$ with $r \geq 1$ and $0 \leq s < q$, and divide the tower $V, TV, \ldots, T^{n-1}V$ into $r$ subtowers of height $q$, except for the topmost $s$ remaining floors, that are left as they are.

Next each subtower is stamped with the stamp we have for $F$ with appropriate width $\mu(V)$, and so we have marked with the stamp some floors along the Rohlin tower, exactly $rp$ of them. All floors that are marked in the Rohlin tower are spaced by runs of $t_j - 1$ unmarked floors, for some $1 \leq j \leq s$. There are $rp_j$ marked floors that are at time $t_j$ to the upper next marked one.

We call $U$ the union of the marked floors in this tower. We have

$$\mu(U) = rp\mu(V) \in \left] \frac{r}{r+1} \frac{p}{q}(1-\delta), \frac{p}{q} \right] \subset \left] (1-\delta)^2 \frac{p}{q}, \frac{p}{q} \right]$$

because $n \geq N$, $(rq+s)\mu(V) \geq 1 - \delta$ and $\frac{1}{\delta} \leq \frac{N}{q}$.

In each subtower that has another one above, say $\widetilde{S}_k = T^{kq}V \cup \cdots \cup T^{kq+q-1}V$, for some $0 \leq k < r$, and for any $t \in ]-\infty, k_s]$, one has

$$\frac{1}{\mu(\widetilde{S}_k)} \mu(x \in \widetilde{S}_k : \alpha \tau_U(x) \leq t) = F(t).$$

So if we denote by $\widetilde{S} = \bigcup_{k<r} \widetilde{S}_k$, we obtain that for any $t \in ]-\infty, k_s]$,

$$\frac{1}{\mu(\widetilde{S})} \mu(x \in \widetilde{S} : \alpha \tau_U \leq t) = F(t).$$

From this the proof follows rather easily because $\mu(U)$ is very close to $\alpha$ and $\mu(\widetilde{S})$ can be made arbitrarily close to 1. $\square$

**5. End of the proof.**

PROOF OF THEOREM 1. Pick an $F \in \mathcal{F}$: by the concavity–continuity lemma, there exists a sequence of rational $F$'s, $(F_n)_{n \geq 1}$, with $F_n \Rightarrow F$.

Then pick a decreasing sequence $\varepsilon_n \downarrow 0$ of positive reals. By the approximation lemma, for each $n \geq 1$, there exists $U_n \subset X$ measurable such that for any $t \geq 0$, there exists an $s \geq 0$, with $|F_n(t) - F_{U_n}(s)| \leq \varepsilon_n$ and $|s - t| \leq \varepsilon_n$. Then obviously $F_{U_n} \Rightarrow F$.

To obtain the reciprocal, assuming $F_{U_n} \Rightarrow F$, the concavity–continuity lemma precisely states that $F$ must belong to $\mathcal{F}$. □

**Acknowledgments.** We would like to thank Sandro Vaienti for proposing the question and motivating preliminary discussions. The authors also thank a first anonymous referee for essential remarks.

Faculty of Mathematics and Physics  
Charles University  
KTIML  
118 00 Prague 1  
Czech Republic  
e-mail: kupsa@kti.mff.cuni.cz

ISITV  
Université du Toulon et du Var  
Avenue Georges Pompidou  
BP 56  
83162 La Valette Cedex  
France  
e-mail: yves.lacroix@univ-tln.fr